\documentclass[11pt]{amsart}
\usepackage{amssymb,latexsym,eufrak,amsmath,amscd,
graphicx}

\theoremstyle{plain}
\newtheorem{theorem}{Theorem}[section]
\newtheorem{proposition}[theorem]{Proposition}
\newtheorem{corollary}[theorem]{Corollary}
\newtheorem{lemma}[theorem]{Lemma}

\theoremstyle{definition}

\newtheorem{remark}[theorem]{Remark}

\newcommand{\PP}{\mathbf{P}}

\newcommand{\NN}{\mathbf{N}}

\newcommand{\QQ}{\mathbf{Q}}

\newcommand{\OO}{\mathcal  {O}}

\newcommand{\JJ}{\mathcal  {J}}

\newcommand{\FF}{\mathcal  {F}}

\newcommand{\fra}{\frak{a}}

\newcommand{\frb}{\frak{b}}

\begin{document}

\title{Global division of cohomology classes via injectivity}

\author[L.Ein]{Lawrence Ein}
\address{Department of Mathematics \\ University of
Illinois at Chicago, \hfil\break\indent  851 South
Morgan Street (M/C 249)\\ Chicago, IL 60607-7045, USA}
\email{ein@math.uic.edu}

\author[M. Popa]{Mihnea Popa}
\address{Department of Mathematics \\ University of
Illinois at Chicago, \hfil\break\indent  851 South
Morgan Street (M/C 249)\\ Chicago, IL 60607-7045, USA}
\email{mpopa@math.uic.edu}

\thanks{LE was partially supported by the NSF grant DMS-0200278.
MP was partially supported by the NSF grant DMS-0500985 and by an AMS Centennial Fellowship.}

\maketitle

$$Dedicated~to~Mel~Hochster~on~his~65th~birthday$$

\section{Introduction}

The aim of this note is to remark that the injectivity theorems of Koll\'ar and Esnault-Viehweg can be used to give a quick algebraic proof of a strengthening (by dropping the positivity hypothesis) of the Skoda-type division theorem for global sections of adjoint line bundles vanishing along suitable multiplier ideal sheaves proved in \cite{el}, and to extend this result to higher cohomology classes as well (cf. Theorem \ref{generalized_skoda}). For global sections, this is a slightly more general statement of the algebraic version of an analytic result of Siu \cite{siu} based on the original Skoda theorem. In \S4 we list a few consequences of this type of result, like the surjectivity of various multiplication or cup-product maps and the corresponding version of the geometric effective Nullstellensatz.

Along the way, in \S3 we write down an injectivity statement for multiplier ideal sheaves (Theorem 
\ref{multiplier_injectivity}), and its implicit torsion-freeness and vanishing consequences (Theorem \ref{torsion}). They are not required in this generality for the main result here (see the paragraph below), but having them available will hopefully be of use. Modulo some standard tricks, the results in \S3 reduce quickly to theorems of Koll\'ar \cite{kollar1} and Esnault-Viehweg \cite{ev}, and we do not claim originality in any of the proofs.

All of the results are proved in the general setting of twists by nef and abundant (or good) line bundles, which replace twists by nef and big line bundles required for the use of vanishing theorems. In particular, what is used in the proof of the main Skoda-type statement is a Koll\'ar vanishing theorem for the higher direct images of adjoint line bundles of the form $K_X + L$, where $L$ is the round-up of a nef and abundant $\QQ$-divisor. For such vanishing, the only contribution we bring here is a natural statement that seems to be slightly more general that what we found in the literature  (cf. Corollary \ref{nef_and_abundant_injectivity}(4)). The proof is otherwise standard, after establishing a simple lemma on restrictions of nef and abundant divisors in \S2.

Mel Hochster has given a beautiful treatment to local statements of Brian\c con-Skoda-type in positive characteristic, using tight closure techniques. We are very happy to be able to contribute work in a similar circle of ideas to a volume in his honor.

\section{Preliminaries}

\noindent
{\bf Injectivity.}
We always work with varieties defined over an algebraically closed field of characteristic zero. 
We first recall that the approach to vanishing theorems described by Esnault and Viehweg in \cite{ev} produces  the following injectivity statement\footnote{Note that \cite{ev} 5.1 contains a slightly more general statement which allows for adding an extra effective divisor under some transversality conditions. We will not make use of this here.}: 

\begin{theorem}[\cite{ev} 5.1]\label{injectivity}
Let $X$ be a smooth projective variety, and let $L$ be a line bundle on $X$. Assume that there exists a
reduced simple normal crossings divisor $\sum_i \Delta_i$ such that we can write 
$L \sim_{\QQ} \sum_i \delta_i  \Delta_i$, with $0 < \delta_i <1$ for all $i$.   
If $B$ is any effective divisor supported on $\sum_i \Delta_i$, then the natural maps 
$$H^i (X, \OO_X(K_X + L)) \longrightarrow H^i (X, \OO_X(K_X + L + B))$$ 
are injective for all $i$. 
\end{theorem}

\noindent
{\bf Nef and abundant divisors.}
Recall next that a nef $\QQ$-divisor $D$ is called \emph{abundant} (or \emph{good}) if $\kappa(D) = \nu (D)$, i.e. its Iitaka dimension is equal to its numerical dimension (the largest integer $k$ such that $D^k\cdot Y \neq 0$ for some $Y \subseteq X$ of dimension $k$).  Note that a semiample divisor is abundant, and so is a big and nef one. Recall the following basic description due to Kawamata:

\begin{lemma}\label{abundant}(\cite{kawamata} Proposition 2.1 and \cite{ev} Lemma 5.11)
For a  $\QQ$-divisor $B$ on a normal projective variety $X$ the following are equivalent:

\noindent
(1) $B$ is nef and abundant. 

\noindent
(2) There exists a birational morphism $\phi: W \rightarrow X$ 
with $W$ smooth and projective, and a fixed effective divisor $F$ on $W$, such that for any $m$ sufficiently large and divisible we can write $\phi^*(mB) = A+ F$, with $A$ semiample. 
\end{lemma}

The following simple restriction statement will be used to note that Koll\'ar vanishing holds for higher 
direct images of $K_X + B$ with $B$ nef and abundant without restrictions, a fact needed here, and slightly stronger than results stated in the literature (cf. Theorem \ref{torsion} (3) and Corollary \ref{nef_and_abundant_injectivity} (4)). 

\begin{lemma}\label{restriction}
Let $X$ be a normal projective variety and $B$ a nef and abundant $\QQ$-divisor on $X$. 

\noindent
(1) If $L$ is a globally generated line bundle on $X$, and $D\in |L|$ a general divisor, then $B_{|D}$ is also nef and abundant.

\noindent
(2) If $C$ is another nef and abundant $\QQ$-divisor on $X$, then $B + C$ is also nef and abundant.

\noindent
(3) If $f:Y \rightarrow X$ is a surjective morphism from another normal projective variety, 
then $f^* B$ is nef and abundant.
\end{lemma}
\begin{proof}
For (1), we only need to show that $B_{|D}$ is abundant. This follows most easily from the description above. Consider a birational morphism $\phi: W \rightarrow X$ as in Lemma \ref{abundant}. If $\tilde{D}$ is the proper transform of $D$, we can choose $D$ such that $\phi_{|D}$ is an isomorphism at the generic point of $D$, $\tilde{D}$ is smooth, and no component of $\tilde{D}$ is contained in ${\rm Supp}(F)$. In this case we have an induced decomposition 
$(\phi_{| \tilde D})^* (mB_{|D}) = A_{| \tilde{D}} + F_{| \tilde{D}}$ for all $m$ sufficiently large and divisible. By Lemma \ref{abundant}, this shows that $B_{|D}$ is abundant. Parts (2) and (3) follow immediately from the same characterization in Lemma \ref{abundant}.
\end{proof}

\begin{remark}
Although we will not use this, it is worth pointing out the following more precise statement, which can be obtained by a closer analysis of the characterization in \cite{kawamata} Proposition 2.1: in the setting of Lemma \ref{restriction}, assume that $\kappa(B) = k$. We always have $B^k \cdot L^{n-k} \ge 0$, and when the restriction $B_{|D}$ is nef and abundant, then $\kappa(B_{|D}) = k$ iff $B^k \cdot L^{n-k} > 0$, while 
$\kappa(B_{|D}) = k-1$ iff $B^k \cdot L^{n-k} = 0$.
\end{remark}

\section{Injectivity for $\QQ$-divisors and multiplier ideals}

In this section we write down the proof of an injectivity statement for multiplier ideals.
This is a consequence of Theorem \ref{injectivity}, and it follows 
quickly from it via standard tricks. For the appropriate level of generality, we use multiplier ideals associated to ideal sheaves. The equivalent divisorial condition (on the log-resolution) is stated at the beginning of the proof. Note that all throughout by 
$\QQ$-effective we mean $\QQ$-linearly equivalent to an effective divisor.

\begin{theorem}\label{multiplier_injectivity}
Let $X$ be a smooth projective variety, $\fra \subseteq \OO_X$ an ideal sheaf, and $L$ a line bundle on $X$. Consider also $A$ a line bundle with $A\otimes \fra$ globally generated, $B$ an effective divisor, and $\lambda\in \QQ$ such that:
\begin{itemize}
\item $L - \lambda A$ is nef and abundant and   
\item $L - \lambda A -\epsilon B$ is $\QQ$-effective, for some $0 < \epsilon < 1$. 
\end{itemize}
Then the natural maps 
$$H^i (X, \OO_X(K_X + L)\otimes \JJ(X, \fra^{\lambda})) \longrightarrow H^i (X, \OO_X(K_X + L + B)
\otimes \JJ(X, \fra^{\lambda}))$$ 
are injective for all $i$. 
\end{theorem}
\begin{proof}
Let $f: Y \rightarrow X$ be a log-resolution of the ideal $\fra$, with $\fra \cdot \OO_Y = \OO_Y(-E)$. The two hypotheses in the Theorem imply that:
\begin{itemize}
\item $f^*L - \lambda E = f^*(L - \lambda A) + \lambda (f^*A - E)$ is nef and abundant 
(by Lemma \ref{restriction}(2) and (3))  
\item $f^*L - \lambda E -\epsilon f^*B$ is $\QQ$-effective, for some $0 < \epsilon < 1$. 
\end{itemize}
We recall that by definition $\JJ(X, \fra^{\lambda}) = f_* \OO_Y (K_{Y/X} - [\lambda E])$.
The projection formula and the Local Vanishing theorem (cf. \cite{positivity} 9.4.4) imply that we have isomorphisms
$$H^i (Y,  \OO_Y ( K_Y + f^*L -[\lambda E])) \cong H^i(X, \OO_X (K_X + L) \otimes \JJ(X, \fra^{\lambda}))$$
and the analogues with $B$ added. Thus it is enough to show injectivity on $Y$, namely for the maps
$$H^i (Y, \OO_Y(K_Y + f^*L -[\lambda E])) \longrightarrow H^i (Y, \OO_Y(K_Y + f^*L - [\lambda E] + f^*B)).$$

By assumption there exists an $a\in \NN$ such that $a(f^*L- \lambda E - \epsilon f^*B) \sim B^{\prime}$, where $B^{\prime}$ is an integral effective divisor. 
On the other hand $f^*L-\lambda E$ is nef and abundant, and we may assume that the log-resolution $f$ factors through the birational morphism $\phi$ of Lemma \ref{abundant}. Consequently we can write 
$f^*L-\lambda E = A^{\prime} + F$, with $A^{\prime}$ a semiample $\QQ$-divisor and $F$ an effective $\QQ$-divisor with fixed support but arbitrarily small coefficients. In particular, for $N>>0$ we can write $(N-a)A^{\prime} \sim P$, where $P$ is a reduced, irreducible divisor.  We can then write down a decomposition
\begin{equation}\label{decomposition}  
f^*L- \lambda E = \frac{1}{N} (a(f^*L- \lambda E) + (N-a) A^{\prime} + (N-a)F) \sim_{\QQ}$$
$$\sim_{\QQ}  \alpha f^*B + \beta B^{\prime} + \gamma P + \frac{N-a}{N} F,
\end{equation}
with $\alpha$, $\beta$ and $\gamma$ arbitrarily small. 

We claim that by passing to a log-resolution of the pair $(X, E + f^*B + B^{\prime}+ P + F)$, we can assume in addition that everything is in simple normal crossings. Let's assume this in order to conclude.
Denote $\Delta' := \lambda E - [\lambda E]$. Using (\ref{decomposition}), we can write:
$$f^*L - [\lambda E] \sim_{\QQ} \Delta^{\prime} + \alpha f^* B + \beta  B^{\prime} + \gamma P 
+ \frac{N-a}{N} F.$$
Note that $\Delta^{\prime}$ may have common components with some of the other divisors appearing 
in the expression on the right hand side. However, since their coefficients can be made arbitrarily small, we can assume that every irreducible divisor in the sum appears with coefficient less than $1$.
Consequently we can apply Theorem \ref{injectivity}, with the role of $L$ played by 
$f^*L - [\lambda E]$, and that of $B$ by $f^*B$. 

It remains to prove the above claim. To this end, note that with our choices we have the following: if we denote
$$T : = \Delta^{\prime} + \alpha f^* B + \beta  B^{\prime} + \gamma P 
+ \frac{N-a}{N} F,$$
then $(Y, T)$ is a klt pair, and we saw that $f^*L - [\lambda E] \sim_{\QQ} T$.
Let $g : Z \rightarrow Y$ be a log-resolution of $(Y, T)$. We have then that 
$$\JJ( Y, T) = g_* \OO_Z (K_{Z/Y} - [g^* T]) \cong \OO_Y.$$ 

Applying again the Local Vanishing theorem cited above, we see then that it is enough to prove injectivity on $Z$, for the map
$$H^i (Z, \OO_Z(K_Z + g^*(f^*L - [\lambda E]) -[g^* T])) \longrightarrow H^i (Z,  \OO_Z(K_Z + g^*(f^*L - [\lambda E])  - [g^* T] + g^*f^*B)).$$ 
But $g^*(f^*L - [\lambda E]) -[g^* T] \sim_{\QQ} \{g^* T\}$, 
which reduces us to the simple normal crossings situation. 
\end{proof}

Theorem \ref{multiplier_injectivity} implies standard torsion-freeness and vanishing consequences for images of twisted multiplier ideal sheaves, in analogy with Koll\'ar's \cite{kollar1} Theorem 2.1.

\begin{theorem}\label{torsion}
Let $X$ be a smooth projective variety, and $f:X \rightarrow Y$ a surjective morphism to a projective variety $Y$. Let $\fra\subset \OO_X$ be an ideal sheaf, and $A$ a line bundle on $X$ such that $A\otimes \fra$ is globally generated. Let $\lambda \in \QQ$ and $L$ a line bundle on $X$ such that $L - \lambda A$ is nef and abundant. Then:

\noindent 
(1) $R^i f_* (\OO_X(K_X + L)\otimes \JJ(\fra^{\lambda}))$ are torsion-free for all $i$.
\newline
\noindent
(2) $R^i f_* (\OO_X(K_X + L)\otimes \JJ(\fra^{\lambda})) = 0 $ for $i > {\rm dim}(X) - {\rm dim}(Y)$.
\newline
\noindent
(3)  $H^j (Y, R^i f_* (\OO_X(K_X + L)\otimes \JJ(\fra^{\lambda}))\otimes M) = 0$ for all $i$ and all $j > 0$, where $M$ is any big and nef line bundle on $Y$. 
\end{theorem}
\begin{proof} 
For the convenience of the reader we sketch briefly how this can be deduced from injectivity, in this case Theorem \ref{multiplier_injectivity}. All the main ideas are of course contained in \cite{kollar1} and \cite{ev}.

The assertion in (2) is an immediate consequence of (1) and base change. For (1), consider $N$ a sufficiently positive line bundle on $Y$. If $R^i f_* (\OO_X(K_X + L)\otimes \JJ(\fra^{\lambda}))$ had torsion, we would be able to choose $D\in |N|$ such that the natural map on global sections 
$$H^0 (Y, R^i f_* (\OO_X(K_X + L)\otimes \JJ(\fra^{\lambda}))\otimes N) \rightarrow 
H^0 (Y, R^i f_* (\OO_X(K_X + L)\otimes \JJ(\fra^{\lambda}))\otimes N(D))$$
is not injective.  On the other hand, for $N$ positive enough, by the degeneration of the corresponding Leray spectral sequences we see that the map above is the same as the natural 
homomorphism
$$H^i (X, \OO_X (K_X + L + f^*N) \otimes \JJ(\fra^{\lambda})) \rightarrow 
H^i (X, \OO_X (K_X + L + f^*N + f^*D) \otimes \JJ(\fra^{\lambda}) ).$$
Note that $f^*N $ is semiample on $X$, and so by Lemma \ref{restriction}(2) the $\QQ$-divisor $L  + f^*N - \lambda A$ is still nef and abundant. We can then apply Theorem \ref{multiplier_injectivity} to derive a contradiction, since the other condition in the theorem is obviously satisfied. 

For (3) one uses induction on the dimension of $X$. We only sketch the proof in the case 
$M$ is ample, which implies the big and nef case via a standard use of Kodaira's Lemma.
For some integer $p >>0$, consider $Y_0 \in |pM|$ a general divisor, so that $X_0 = f^{-1} (Y_0)$ is a smooth divisor in $X$. Fix an $i$, and for simplicity let us denote $\FF : = \OO_X(K_X + L)\otimes \JJ(\fra^{\lambda})$. We have an exact sequence
$$0 \longrightarrow \FF \longrightarrow \FF(X_0) \longrightarrow \FF(X_0)_{|X_0} \longrightarrow 0.$$
Pushing this forward to $Y$, for each $i$ we get exact sequences
$$0 \longrightarrow R^i f_* \FF \longrightarrow R^i f_* (\FF(X_0)) \longrightarrow R^i f_* (\FF(X_0)_{|X_0}) \longrightarrow 0.$$
The reason these sequences are exact at the extremities is that by (1) the sheaves $R^i f_* \FF$ are 
torsion-free, while $R^i f_* (\FF(X_0)_{|X_0})$ are generically zero. Twisting this with $M$, and recalling that $p>>0$, we immediately obtain
$$H^{j+1} (Y, R^i f_* \FF \otimes M) \cong H^j (Y_0, R^i f_* (\FF(X_0)_{|X_0})  \otimes M_{Y_0}), {\rm ~for~all~} j \ge 1.$$
Since $Y_0$ was chosen general, the restriction $(L- \lambda A)_{|X_0}$ is still nef and abundant
by Lemma \ref{restriction}(1), while by \cite{positivity} 9.5.35 we have that $\JJ(\fra^{\lambda})\cdot \OO_{X_0} \cong \JJ(\fra^{\lambda}\cdot \OO_{X_0})$. Thus, by induction on the dimension, we can apply (3) to the right hand side above, and deduce the conclusion on $Y$ for $j \ge 2$. 

We are left with the case $j = 1$. Consider again $p >>0$, and choose a divisor $D \in |p M|$, so that 
$H^1 (Y, R^i f_* \FF \otimes M (D)) = 0$. We have a Leray spectral sequence 
$$E^{j,i}_2:= H^j (Y, R^i f_* \FF \otimes M) \Rightarrow H^{i +j} (X, \FF\otimes f^* M).$$ 
By the previous paragraph the $E^{j,i}_2$ are zero for $j \ge 2$, which by chasing the spectral sequence easily gives that we have an injective map 
$$H^1 (Y, R^i f_* \FF \otimes M) \longrightarrow H^{i +1} (X, \FF\otimes f^* M)$$
Analogously, there is a similar injective map after twisting with $D$. But Theorem \ref{multiplier_injectivity} provides the injectivity of the map
$$H^{i +1} (X, \FF\otimes f^* M) \longrightarrow H^{i +1} (X, \FF\otimes f^* M(D)),$$ 
so finally all of this implies that the map
$$H^1 (Y, R^i f_* \FF \otimes M) \longrightarrow H^1 (Y, R^i f_* \FF \otimes M(D))$$ 
is also injective. However this last group is zero.
\end{proof}

A special case of the two theorems in this section is the following result for $\QQ$-divisors. The first part is  already explicitly stated by Esnault-Viehweg \cite{ev}, and is an extension of Koll\' ar's Injectivity Theorem \cite{kollar1}, Theorem 2.2. Parts (2) and (3) of course follow directly from it. Part (4) is stated a little less generally in \cite{ev} 6.17(b) -- we will need however the fomulation below.

\begin{corollary}\label{nef_and_abundant_injectivity} 
Let $L$ a line bundle on a smooth projective $X$, and $\Delta = \sum_i \delta_i \Delta_i$ a simple normal crossings divisor with $0 < \delta_i <1$ for all $i$. Assume that $L -\Delta$ is nef and abundant, and $B$ is an effective Cartier divisor such that $L - \Delta - \epsilon B$ is $\QQ$-effective for some $0 < \epsilon < 1$. Consider also a morphism $f:X\rightarrow Y$, with $Y$ projective. Then:

\noindent
(1) (\cite{ev} 5.12(b).) The natural maps 
$$H^i (X, \OO_X(K_X + L)) \longrightarrow H^i (X, \OO_X(K_X + L + B))$$ 
are injective for all $i$. 

\noindent
(2) $R^i f_* \OO_X(K_X + L)$ are torsion-free for all $i$.

\noindent
(3) $R^i f_* \OO_X(K_X + L) = 0 $ for $i > {\rm dim}(X) - {\rm dim}(Y)$.

\noindent
(4)  $H^j (Y, R^i f_* \OO_X(K_X + L)\otimes M) = 0$ for all $i$ and all $j > 0$, where $M$ is any big and nef line bundle on $Y$. 
\end{corollary}

\section{Skoda-type global division theorem}

For the objects involved in the statement below recall the following. If $f:Y \rightarrow X$ is a common log-resolution for ideal sheaves $\fra$ and $\frb$, with $\fra \cdot \OO_Y = \OO_Y(-E)$ and $\frb\cdot \OO_Y = \OO_Y(-F)$, and if $\mu, \lambda \in \QQ_{+}$, then the ``mixed" multiplier ideal is defined as $\JJ(\fra^{\mu} \cdot \frb^{\lambda}) = f_* \OO_Y(K_{Y/X} - [\mu E + \lambda F])$ (cf. \cite{positivity} 9.2.8). For a line bundle $B$ on $X$, we say that $B\otimes \frb^{\lambda}$ is \emph{nef and abundant} if the $\QQ$-divisor $f^*B - \lambda F$ is nef and abundant on $Y$. These definitions are independent of the log-resolution we choose.

\begin{theorem}\label{generalized_skoda}
Let $X$ be a smooth projective variety of dimension $n$, and $\fra, \frb \subseteq \OO_X$ ideal sheaves. Consider line bundles $L$ and $B$ on $X$ such that $L\otimes \fra$ is globally generated, and $B\otimes \frb^{\lambda}$ is nef and abundant for some $\lambda \in \QQ_{+}$. Then for every integer $m\ge n+2$, the sections in 
$$H^0 (X, \OO_X (K_X + mL + B) \otimes \JJ(\fra^m \cdot \frb^{\lambda}))$$  
can be written as linear combinations (with coefficients in $H^0(L)$) of sections in 
$H^0(X, \OO_X (K_X + (m-1)L + B) \otimes \JJ(\fra^{m-1} \cdot \frb^{\lambda}))$.
More generally, for every $i\ge 0$, the cohomology classes in 
$$H^i (X, \OO_X (K_X + mL + B) \otimes \JJ(\fra^m \cdot \frb^{\lambda}))$$  
can be written as linear combinations (with coefficients in $H^0(L)$, via cup product) of classes in 
$H^i(X, \OO_X (K_X + (m-1)L + B) \otimes \JJ(\fra^{m-1} \cdot \frb^{\lambda}))$.
\end{theorem}
\begin{proof}
Let $f: Y \rightarrow X$ be a common log-resolution for $\fra$ and $\frb$, with $\fra\cdot \OO_Y = \OO_Y(-E)$, and $\frb\cdot \OO_Y = \OO_Y(-F)$. Let's assume that $L\otimes \fra$ is generated by 
$s \ge n+1$ sections spanning a linear subspace $V \subseteq H^0L$. In fact, it will be clear from the proof that if $L\otimes \fra$ happens to be spanned by $p \le n$ sections, then the result can be improved, with an identical argument, replacing $n+1$ by $p$. 

The line bundle $A: = f^* L - E$ will be generated by the corresponding space of sections in $H^0(A)$, which we denote also by $V$.  
Denote by $g: Y \rightarrow Z\subset \PP^{s-1}$ the map it determines, so that $A \cong g^* \OO_Z(1)$. 
The goal is to prove the surjectivity of the multiplication map:
$$V \otimes H^0 (X, \OO_X(K_X + (m-1)L + B) \otimes \JJ(\fra^{m-1} \cdot \frb^{\lambda}))
\rightarrow H^0 (X, \OO_X (K_X + mL + B) \otimes \JJ(\fra^m \cdot \frb^{\lambda})).$$ 
Note however that by definition 
$$\JJ(\fra^k \cdot \frb^{\lambda}) = f_* \OO_Y(K_{Y/X} - kE - [\lambda F]) {\rm ~for~all ~}k,$$ 
and so  this is equivalent to the surjectivity of the multiplication map
$$V \otimes H^0 (Y, \OO_Y(K_Y + (m-1)f^*L + f^*B - (m-1)E - [\lambda F]))
\rightarrow H^0 (Y, \OO_Y (K_Y + mf^*L + f^*B - mE - [\lambda F])).$$
Recall the notation $A = f^*L - E$, and denote $N = f^*B - [\lambda F]$, which by assumption can be written as a nef and abundant $\QQ$-divisor plus a simple normal crossings boundary divisor. 
Rewritting the above, we are then interested in the surjectivity of the multiplication map
\begin{equation}\label{goal}
V \otimes H^0 (Y, \OO_Y(K_Y + (m-1)A + N)) \rightarrow  H^0 (Y, \OO_Y(K_Y + mA + N)).
\end{equation}

We compare this with the picture obtained by pushing forward to $Z$ via $g$. Note that $V$ can be considered as a space of sections generating 
$\OO_Z(1)$. It is well known that this gives rise to an exact Koszul complex on $Z$ (cf. e.g. \cite{positivity}, beginning of Appendix B.2):
$$ 0 \rightarrow \Lambda^{s} V \otimes \OO_Z(-s) \rightarrow \ldots \rightarrow \Lambda^2 V \otimes \OO_Z(-2) \rightarrow V\otimes \OO_Z(-1) \rightarrow \OO_Z \rightarrow 0.$$
We twist this with the sheaf $g_*\OO_Y(K_Y + N) \otimes \OO_Z(m)$, which preserves the exactness of the sequence:\footnote{Note that the Koszul complex is locally split, and its syzygies are locally free, 
so twisting by any coherent sheaf preserves exactness.}
$$0 \rightarrow \Lambda^{s} V \otimes \OO_Z(m-s)\otimes g_*\OO_Y(K_Y + N) \rightarrow \ldots \rightarrow $$
$$\rightarrow V\otimes \OO_Z(m-1) \otimes g_*\OO_Y(K_Y + N ) \rightarrow \OO_Z(m) \otimes g_* \OO_Y (K_Y + N) \rightarrow 0 .$$
Since $A \cong g^* \OO_Z(1)$, the surjectivity of the map in $(\ref{goal})$ is equivalent to the surjectivity of the multiplication map
$$V\otimes H^0(Z, \OO_Z(m-1) \otimes g_*\OO_Y(K_Y + N )) \rightarrow H^0 (Z, \OO_Z(m) \otimes g_* \OO_Y (K_Y + N)) $$
induced by the taking global sections in the Koszul complex above. 
But since $m \ge n+2$, by the line bundle case of Corollary \ref{nef_and_abundant_injectivity} (4) we have 
$$H^j (Z, \OO_Z (m- i) \otimes g_* \OO_Y (K_Y + N)) = 0, ~{\rm for~all~} j>0 
{\rm~and~} i \le n+1.$$
By chasing through the induced short exact sequences, this easily implies that the entire Koszul complex stays exact after passing to global sections. (Note that beyond $i = n+1$ we are only interested in cohomology groups as above for $j >n$, which are automatically $0$.) This proves the statement for global sections. 

The proof of the general statement for cohomology classes is similar. Note first that again because of Corollary \ref{nef_and_abundant_injectivity} (4), for higher direct images, we have that for every $i$ and every $k>0$, the Leray spectral sequence degenerates to an isomorphism
$$H^i (Y, \OO_Y(K_Y + kA + N)) \cong H^0 (Z, R^i g_* \OO_Y (K_Y + N) \otimes \OO_Z(k)).$$
But, exactly as above, the same vanishing applied yet again, for $R^i g_*\OO_Y(K_Y + N)$, implies that after twisting the Koszul complex with
$R^i g_* \OO_Y (K_Y + N)$ and passing to global sections, we obtain a surjection
$$V\otimes H^0(Z, \OO_Z(m-1) \otimes R^i g_*\OO_Y(K_Y + N )) \rightarrow H^0 (Z, \OO_Z(m) \otimes R^i g_* \OO_Y (K_Y + N)).$$
This implies that we have the surjectivity of the cup product
$$V \otimes H^i (Y, \OO_Y(K_Y + (m-1)A + N)) \rightarrow  H^i (Y, \OO_Y(K_Y + mA + N)).$$
On the other hand, by the Local Vanishing theorem (cf. \cite{positivity} 9.4.4) we have that 
$$R^j f_* \OO_Y ( K_Y - k E - [\lambda F]) = 0, {\rm ~for~all~} j > 0 {\rm ~and~} k>0$$ 
so consequently, for all $i$,
$$H^i (X, \OO_X (K_X + kL + B) \otimes \JJ(\fra^m \cdot \frb^{\lambda})) \cong 
H^i (Y, \OO_Y(K_Y + kA + N))$$ 
and the result follows.
\end{proof}

By taking $\frb = \OO_X$ in the Theorem we obtain the following: 

\begin{corollary}\label{siu}
In the notation of Theorem \ref{generalized_skoda}, if $B$ is a nef and abundant line bundle (for example $B = \OO_X$), we have that for every $m\ge n+2$ and every $i$, the cohomology classes in 
$$H^i (X, \OO_X (K_X + mL + B) \otimes \JJ(\fra^m))$$  
can be written as linear combinations of classes in 
$H^i(X, \OO_X (K_X + (m-1)L+ B) \otimes \JJ(\fra^{m-1}))$. 
\end{corollary}
 
In the case of global sections, i.e. $i = 0$, this is an improvement of the Global Divison Theorem  of \cite{el} (cf. \cite{positivity}, Theorem 9.6.31): one does not require twisting with an ample (or big and nef) line bundle.\footnote{Note though that in order to do this we have to start with $m = n+2$ as opposed to $m= n+1$.} Also for $i = 0$, it is a slightly more general version of the algebraic version of the Skoda-type theorem proved in the anaytic context in \cite{siu} Theorem 1.8.3. 
Note that the method used here does not distinguish between global sections and higher cohomology classes. The case of trivial multiplier ideals already yields a slightly surprising statement even in the case $B = \OO_X$.

\begin{corollary}\label{free}
If $L$ is a globally generated line bundle and $B$ is a nef and abundant line bundle, on a smooth projective variety $X$ of dimension $n$, then for all $m \ge n+2$ and all $i \ge 0$ the cup product maps 
$$H^0(X, L)\otimes H^i(X, \OO_X(K_X + (m-1)L + B))\longrightarrow H^i(X, \OO_X(K_X + mL + B))$$
are surjective.
\end{corollary}

An interesting consequence of this involves multiplication maps of globally generated adjoint line bundles. When $B$ is globally generated, this is the statement of \cite{siu} Theorem 1.8.4. 

\begin{corollary}
Let $B$ be a nef and abundant line bundle on $X$ such that the adjoint bundle $L:= K_X + B$ is globally generated. Then $mL$ is projectively normal for all $m\ge n+2$ and the section ring 
$R_L = \underset{m\ge 0}{\bigoplus} H^0 (mL)$ is generated by $\underset{m\le n+2}{\bigoplus} H^0 (mL)$.\footnote{In particular, if $L$ is also ample, then $mL$ is very ample for $m\ge n+2$. This can be shown directly though using Castelnuovo-Mumford regularity.}
\end{corollary}
\begin{proof}
Corollary \ref{free} implies that for all $m \ge n+2$ we have the surjectivity the multiplication map
$$H^0 (L) \otimes H^0 (mL) \longrightarrow H^0 ((m+1) L),$$
which implies the generation statement. By iteration we obtain the surjectivity of the map
$$H^0 (L)^{\otimes k} \otimes H^0 (mL) \longrightarrow H^0 ((m+k) L),$$
for all $k \ge 1$, which has as a special consequence the projective normality of  $mL$.
\end{proof}

For completeness, note that the statement of Corollary \ref{free}, at least for $B = \OO_X$, holds also for higher direct images  of canonical bundles.

\begin{proposition}\label{higher}
Let $X$ and $Y$ be projective varieties, with $X$ smooth and $Y$ of dimension $n$, and let $f: X \rightarrow Y$ be a surjective morphism. Consider a globally generated line bundle $L$ on $Y$. 
Then, for all $m \ge n+2$ and all $i, j \ge 0$, the cup product maps 
$$H^0(Y, L)\otimes H^i(Y, R^j f_* \omega_X \otimes \OO_Y((m-1)L))\longrightarrow H^i(Y, R^j f_* \omega_X \otimes \OO_Y(mL))$$
are surjective.
\end{proposition}
\begin{proof}
The proof goes along the same lines as the proof of Theorem \ref{generalized_skoda}. We use the 
morphism $g: Y \rightarrow Z$ induced by $L$ and the corresponding Koszul complex on $Z$. The 
extra thing to note is that the vanishing
$$H^k (Z, R^i g_* R^j f_* \omega_X \otimes \OO_Z (l)) = 0, \forall k, l >0$$
still holds, according to \cite{kollar2} Theorem 3.4.
\end{proof}

One also obtains a similar weakening of the hypotheses under which the Geometric Effective Nullstellensatz, \cite{el} Theorem (iii) (cf. also \cite{positivity} Theorem 10.5.8) holds, by plugging the statement of Theorem \ref{generalized_skoda} into the original proof of \emph{loc. cit.}

\begin{corollary}[Geometric effective Nullstellensatz]
Let $X$ be a smooth projective variety, $\fra \subset \OO_X$ an ideal sheaf, and $L$ an ample line bundle on $X$ such that $L \otimes \fra$ is globally generated, say by sections $g_j \in H^0 (L\otimes 
\fra)$. Consider also a nef and abundant line bundle $B$ on $X$. Then, for all $m \ge n+2$, if a section 
$$g \in H^0 (X, \OO_X(K_X + mL + B)) $$
vanishes to order at least $(n+1)\cdot {\rm deg}_L X$ at a general point of each distinguished subvariety\footnote{Recall that this means the following: if $f:Y \rightarrow X$ is the normalized blow-up of $X$ along $\fra$, and $\fra \cdot \OO_Y = \OO_Y ( -\sum a_i E_i)$, then the distinguished subvarieties of $\fra$ are the images of the $E_i$'s in $X$ -- cf. \cite{el} \S2.} of $\fra$, then 
$g$ can be expressed as a linear combination $\sum h_j g_j$, with 
$h_j \in H^0 (X, \OO_X(K_X + (m-1)L + B))$. 
\end{corollary}

\noindent
{\bf Acknowledgements.} We would like to thank the referee for a very detailed reading of the paper 
and for useful suggestions.

\end{document}